\documentstyle[amscd, amstex]{amsart}

\def\ZZ         {{\Bbb Z}}
\def\RR         {{\Bbb R}}
\def\CC         {{\Bbb C}}

\def\PP         {{\Bbb P}}

\def\AA         {{\Bbb A}}

\newtheorem{thm}{Theorem}[section]

\newtheorem{pr}[thm]{Proposition}

\theoremstyle{definition}
\newtheorem{rem}[thm]{Remark}
\newtheorem{ex}[thm]{Example}
\newtheorem{defn}[thm]{Definition}

\newcommand{\mult}{{\rm mult}}

\newcommand{\spe}{{\rm Spec}}
\newcommand{\key}{\bibitem}

\newcommand{\ps}{{{\bf P}_{\Sigma}}}

\newcommand\hidot{{\raise1pt\hbox{$\scriptscriptstyle\bullet$}}}
\newcommand\lodot{{\raise.3pt\hbox{$\scriptscriptstyle\bullet$}}}
\newcommand{\dd}{{\rm d}}
\newcommand{\pp}{{\bf P}}

\begin{document}

\title[Deformations of Calabi-Yau hypersurfaces]
{Deformations of Calabi-Yau hypersurfaces arising from
deformations of toric varieties}
\author{Anvar R. Mavlyutov}
\address {Department of Mathematics, Indiana University, Bloomington, IN 47405, USA.}
 \email{amavlyut@@indiana.edu}

\keywords{Deformations, Calabi-Yau, toric varieties, mirror
symmetry.}
\subjclass{Primary: 14M25}

\maketitle

\tableofcontents

\setcounter{section}{-1}

\section{Introduction.}

It is a well known fact now  that the infinitesimal deformations
of Calabi-Yau manifolds are unobstructed (see \cite{ti,to}). A
generalization of this theorem for Calabi-Yau orbifolds
nonsingular in codimension two is due to Z.~Ran (see \cite{r}).
There are easy ``polynomial" deformations of Calabi-Yau
hypersurfaces in toric varieties performed by changing the
coefficients of the defining polynomial of the hypersurface. But
according to the unobstructedness and the calculation of the space
of infinitesimal deformations (e.g., \cite{m2}), there must also
exist ``non-polynomial" deformations of the minimal Calabi-Yau
hypersurfaces leading outside the ambient toric variety. In this
paper, we have constructed  the missing non-polynomial
deformations, which are explicitly described as abstract
Calabi-Yau varieties.

The first examples of non-polynomial deformations mentioned in
\cite{cdfkm} were constructed by Sheldon Katz and David Morrison
for two parameter models (e.g., for crepant resolutions of degree
8 hypersurfaces in $\PP(1,1,2,2,2)$). Despite  the importance of
these deformations in string theory (see \cite{kmp}, \cite{kklm}
), a general solution was elusive for quite a while. All of the
previous constructions of non-polynomial deformations were
realized by embedding an ample Calabi-Yau hypersurface of some
weighted projective space into another variety where a deformation
was performed, and an appropriate blow-up of this deformation
produced the deformation of the crepant resolution. Our approach
is different. First, we found a credible evidence that the
non-polynomial deformations of the Calabi-Yau hypersurfaces should
be induced by the deformations of the ambient toric varieties.
Second, the nature of the definition of a toric variety led us to
 the idea to reglue the affine toric varieties in order to
deform the complex structure on the ambient variety. A further
investigation revealed us that this deformation is obtained by an
automorphism of an open toric subvariety corresponding to its
root.  Finally, the deformation gives a flat family of
hypersurfaces changing the complex structure on the original
Calabi-Yau hypersurface.

We have the following plan for the paper. Section 1 gives a set up
of the toric varieties in terms  of homogeneous coordinates, and
also describes the space of infinitesimal deformations of minimal
Calabi-Yau hypersurfaces. Then, Section 2 has the construction of
deformations of complete toric varieties with a big and nef
divisor. We show that such deformations induce deformations of the
hypersurfaces as well. In Section 3, we calculate the
infinitesimal deformations in terms of \v Cech cocycles
corresponding to the global deformations of the hypersurfaces. For
a basis of the space of infinitesimal deformations of the minimal
Calabi-Yau hypersurface, we get the corresponding global
deformations.

{\it Acknowledgement.}  The author would like to thank Sheldon
Katz for explaining the example of non-polynomial deformations of
Calabi-Yau hypersurfaces coming from $\PP(1,1,2,2,2)$.

\section{Infinitesimal deformations of  Calabi-Yau hypersurfaces.}\label{s:i}

This section  recalls  the  definition of a  toric variety as a
quotient of a Zariski open subset of an affine space by a linear
diagonal action of some algebraic subgroup of a complex torus. We
also review the description of the space of infinitesimal
deformations $H^1(X,{\cal T}_X)$ for a minimal Calabi-Yau
hypersurface $X$  in a complete simplicial toric variety $\ps$.
Then, we show that the non-polynomial infinitesimal deformations
represented by \v Cech cocycles  have a ``lift'' to the space
$H^1(\ps,{\cal T}_\ps)$. This leads us to the crucial observation
that the non-polynomial deformations of a Calabi-Yau hypersurface
should be induced by the deformations of the ambient toric
variety. Basic facts about toric varieties can be found in
\cite{c2,f1,d,o}.

We start defining toric varieties with a lattice  $N$
 of rank $d$. A rational convex polyhedral cone $\sigma\subset
N_\RR:=N\otimes_\ZZ \RR$ is a cone generated by finitely many
elements $u_1,\dots,u_s\in N$:
 $$\sigma=\{\lambda_1 u_1+\cdots+\lambda_s
 u_s:\quad\lambda_1\dots\lambda_s\ge0\}.$$
 Such a cone $\sigma$ is strongly convex if
 $\sigma\cap(-\sigma)=0$.
A face  of  $\sigma$ is the intersection of $\sigma$ with one of
its supporting hyperplanes: $\{L=0\}\cap\sigma$, where $L$ is a
linear form, nonnegative on $\sigma$. A {\it fan} $\Sigma$ is
defined to be a finite collection of strongly convex rational
polyhedral cones in $N_\RR$ such that\\
  (i) Each face of a cone in $\Sigma$
belongs to $\Sigma$,\\ (ii) The intersection of two cones in
$\Sigma$ is a face of each.

Denote by $\Sigma(k)$ the set of $k$-dimensional cones of a fan
$\Sigma$. Suppose that 1-dimensional cones of $\Sigma$ span
$N_\RR$. Then, by \cite{c1}, the toric variety $\ps$, associated
with the fan $\Sigma$, can be constructed as follows. Let
$\AA^n={\rm Spec}(\CC[x_1,\dots,x_n])$, where $n$ is the number of
1-dimensional cones $\Sigma(1)$, and let $B(\Sigma)=\langle
\prod_{\rho_i\not\in\sigma} x_i\rangle$ be the ideal in
$\CC[x_1,\dots,x_n]$. Then, the complement to the closed subset
defined by $B(\Sigma)$ gives the Zariski open set
$\AA^n\setminus{\bf V}(B(\Sigma))$. This set is invariant under
the action of an affine algebraic D-group $G$, a subgroup of
$(\CC^*)^n$, defined as a kernel of the group homomorphism
 $(\CC^*)^n@>>>(\CC^*)^d$ sending $(x_1,\dots,x_n)$ to
$(\prod_{i=1}^n x_i^{\langle m_1,e_i\rangle},\dots,\prod_{i=1}^n
x_i^{\langle m_d,e_i\rangle})$, where $e_1,\dots,e_n$ are the
minimal lattice generators of the 1-dimensional cones
$\rho_1,\dots,\rho_n$ of the fan $\Sigma$, and where
$m_1,\dots,m_d$ is a basis of the dual lattice $M={\rm
Hom}(N,\ZZ)$ and $\langle\,,\,\rangle$ is the pairing. Theorem~2.1
in \cite{c1}, allows us to define the toric variety $\ps$,
associated with the fan $\Sigma$, as the categorical quotient of
$\AA^n\setminus{\bf V}(B(\Sigma))$ by $G$. The ring $S={\Bbb
C}[x_1,\dots,x_n]$ is called  the homogeneous coordinate ring
 of $\ps$, where the variables
$x_1,\dots,x_n$ give the irreducible torus $(\CC^*)^n/G$-invariant
divisors $D_1,\dots,D_n$. This ring is graded by the Chow group
$A_{d-1}(\ps)$, assigning $[\sum_{i=1}^n a_i D_i]$ to
$\deg(\prod_{i=1}^n x_i^{a_i})$.

To describe the infinitesimal deformations we first consider a
more general situation than the one of a Calabi-Yau hypersurface.
As in \cite[Section 4]{m2}, let $X$ be a big and nef hypersurface,
defined by $f\in S_\beta$, in the toric variety $\ps$. Then, by
Proposition 1.2 in \cite{m1}, there is the associated toric
morphism $\pi:{\bf P}_\Sigma\rightarrow{\bf P}_{\Sigma_X}$.
Consider a 2-dimensional cone $\sigma\in\Sigma_X$ with at least
one 1-dimensional cone $\rho\subset\sigma$, whose generator lies
in the relative interior ${\rm int}(\sigma)$ of $\sigma$. Using
such a cone $\sigma$ form an open covering of the toric variety
$\ps$ by the sets
 $$U_{\sigma'}=\biggl\{x\in\ps:\prod_{\rho_i\subset\sigma\setminus\sigma'}x_i\ne0\biggr\}$$
 for
all cones $\sigma'\in\Sigma(2)$ that lie in $\sigma$. Fix an order
for this open covering corresponding to as the cones lie inside
$\sigma$:
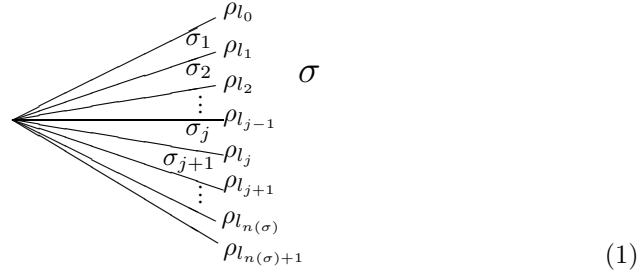
\begin{equation}
\setlength{\unitlength}{1cm}
\begin{picture}(8,3.5)
\put(4.8,3.3){$\rho_{l_0}$} \put(2,1.9){\line(2,1){2.7}}
\put(4.3,2.9){$\sigma_1$} \put(4.8,2.8){$\rho_{l_1}$}
\put(2,1.9){\line(3,1){2.7}} \put(4.3,2.5){$\sigma_2$}
\put(2,1.9){\line(6,1){2.7}} \put(4.8,2.34){$\rho_{l_2}$}
\multiput(4.5,2.2)(0,-.1){3}{\circle*{0.01}}
\put(4.8,1.88){$\rho_{l_{j-1}}$} \put(4.3,1.7){$\sigma_{j}$}
\put(2,1.9){\line(1,0){2.8}} \put(4.8,1.4){$\rho_{l_j}$}
\put(2,1.9){\line(6,-1){2.8}} \put(4,1.3){$\sigma_{j+1}$}
\put(2,1.9){\line(3,-1){2.8}} \put(4.8,1){$\rho_{l_{j+1}}$}
\put(5.8,2.4){\LARGE$\sigma$}
\multiput(4.5,1)(0,-.1){3}{\circle*{0.01}}
\put(4.78,0.52){$\rho_{l_{n(\sigma)}}$}
\put(2,1.9){\line(2,-1){2.7}} \put(2,1.9){\line(5,-3){2.73}}
\put(4.8,0.1){$\rho_{l_{n(\sigma)+1}}$}
\end{picture}\label{e:pic}
\end{equation}
where $n(\sigma)$ is the number of cones $\rho$ such that
$\rho\subset\sigma$ and $\rho\notin \Sigma_X(1)$.

Now let $\ps$ be a simplicial toric variety and assume in addition
that the hypersurface $X\subset\ps$ is quasismooth, i.e., the
partial derivatives  $\partial f/\partial x_{1},\dots,\partial
f/\partial x_{n}$ do not vanish simultaneously on $\ps$. Consider
a refinement $U_{i,\sigma_j}=U_i\cap U_{\sigma_j}$ of the above
open covering and the open covering ${\cal U}=\{U_i\}_{i=1}^n$,
 where $U_i=\{x\in\pp: f_i(x)\ne0\}$ and $f_i:=\partial f/\partial x_{i}$.
 Denote the refined covering
${\cal U}^\sigma$, considering the order on this covering as the
lexicographic order for the pairs of indices $({i,j})$.

In \cite{m2}, we had two different types of \v Cech cocycles which
represent some elements in $H^1(X,{\cal T}_X)$:
\begin{defn} \label{d:non}
Given $A\in S_{\beta}$, set $$(\gamma_A)_{i_0i_1}=
{\fracwithdelims\{\}{\sqrt{-1}A\langle\partial_{i_0}\wedge
\partial_{i_1},\dd f\rangle} {f_{i_0} f_{i_1}}}_{i_0i_1},$$
where $\partial_{i}:=\frac{\partial}{\partial x_{i}}$ and
$\langle\,,\,\rangle$ denotes the contraction.

Given $\rho_l\subset\sigma\in\Sigma_X(2)$ such that
$\rho_l\notin\Sigma_X(1)$, then, as in (\ref{e:pic}), $l=l_j$ for
some $j$, and we set
$$\partial^l_{j}=\frac{x_{l_{j-1}}\partial_{l_{j-1}}}{\mult(\sigma_j)},\quad
\partial^i_{j+1}=-\frac{x_{l_{j+1}}\partial_{l_{j+1}}}{\mult(\sigma_{j+1})},
\text{ and }\partial^l_{k}=0\text{  for }k\ne j,j+1.$$ For $B\in
S_{\beta^\sigma_1}$ (here,
$\beta^\sigma_1:=\sum_{\rho_i\subset\sigma}\deg(x_i)$), define $$
(\gamma^l_B)_{(i_0,{j_0})(i_1,{j_1})}={\biggl\{\frac{B}{\prod_{\rho_i\subset\sigma}
x_i} \biggl(\frac{\langle\partial_{i_1}\wedge\partial^l_{j_1},\dd
f\rangle}{f_{i_1}}-
\frac{\langle\partial_{i_0}\wedge\partial^l_{j_0},\dd
f\rangle}{f_{i_0}}\biggr)\biggr\}} _{(i_0,{j_0}),(i_1,{j_1})}. $$
\end{defn}

The above cocycles give  maps: $\gamma_{\_}:S_\beta\rightarrow
H^1(X,{\cal T}_X)$ and
$\gamma^l_{\_}:S_{\beta^\sigma_1}\rightarrow H^1(X,{\cal T}_X)$.
The following statement is a part of Theorem~7.1 in \cite{m2}.

\begin{pr}\label{t:chiral}
Let $X\subset\ps$ be  a  semiample anticanonical nondegenerate
(i.e., transversal to the torus orbits) hypersurface in a complete
simplicial toric variety, defined by $f\in S_\beta$. Then there is
an isomorphism $$\gamma_{\_}\oplus(\oplus_l\gamma^l_{\_}):
R_1(f)_{\beta} \oplus\Biggl(\bigoplus\begin{Sb}\sigma
\in\Sigma_X(2)\\ e_l\in{{\rm int}(\sigma)}\end{Sb}(S/\langle
x_l\rangle)_{\beta_1^\sigma}\Biggr)\cong H^1(X,{\cal T}_X),$$
where the sum $\oplus_l\gamma^l_{\_}$ is over
$\rho_l\subset\sigma\in\Sigma_X(2)$  such that
$\rho_l\notin\Sigma_X$, and $R_1(f)=S/J_1(f)$ with
$J_1(f):=\langle x_1(\partial f/\partial x_1),\dots,x_n(\partial
f/\partial x_n)\rangle:x_1\cdots x_n$.
\end{pr}

\begin{rem} In the above proposition, if the Calabi-Yau
hypersurface $X$ in $\ps$ has terminal singularities, as in the
case of a maximal projective subdivision of the fan of a Fano
toric variety (see \cite{b}), then $X$ is a minimal Calabi-Yau
orbifold (see Definition 1.4.1 in \cite{ck}), and, by Proposition
A.4.2 in \cite{ck}, the space $H^1(X,{\cal T}_X)$ classifies
infinitesimal deformations.
\end{rem}

The space $R_1(f)_{\beta}$ in Proposition~\ref{t:chiral}
represents the polynomial deformations of a semiample
nondegenerate Calabi-Yau hypersurface performed inside the toric
variety.
 Indeed, the space
$H^0(X,{\cal N}_{X/\ps})$ classifies the infinitesimal
deformations  of $X$ as a subvariety of $\ps$ (see
\cite[Chapter~III,~Exercise~9.7]{h}). On the other hand, the short
exact sequence
 $$0@>>>{\cal T}_X@>>>i^*{\cal T}_\ps@>>>{\cal N}_{X/\ps}@>>>0$$
(here, $i:X\subset\ps$ is the inclusion)        gives the maps
$$H^0(X,{\cal N}_{X/\ps})@>>>H^1(X,{\cal T}_X)@>>>H^1(X,i^*{\cal
T}_\ps).$$
 As  remarked after Definition~3.1 in \cite[Section~3]{m2},
 the restriction of the cocycle  $(\gamma_A)_{i_0i_1}$ to
  $C^p({\cal U}|_X,i^*\wedge^p{\cal T}_\pp)$ is a \v Cech coboundary.
  By the exactness of the above sequence at  the middle term, we deduce that $\gamma_A$ represent the polynomial
deformations, i.e., those that can be performed inside the ambient
space. We will now show that the other part --- the non-polynomial
infinitesimal deformations --- in $H^1(X,{\cal T}_X)$ have a
``lift'' to $H^1(\ps,{\cal T}_\ps)$. The images of $\gamma^l_B$,
for $B\in S_{\beta^\sigma_1}$, in $H^1(X,i^*{\cal T}_\ps)$ are
represented by
 $${\biggl\{\frac{B}{\prod_{\rho_i\subset\sigma}
x_i} \biggl(\partial^l_{j_1}-\partial^l_{j_0}
-\frac{\langle\partial^l_{j_1},\dd
f\rangle}{f_{i_1}}\partial_{i_1}+
\frac{\langle\partial^l_{j_0},\dd
f\rangle}{f_{i_0}}\partial_{i_0}\biggr)\biggr\}}
_{(i_0,{j_0})(i_1,{j_1})}. $$
 Since the polynomials $\langle\partial^l_{j_0},\dd
f\rangle$ are divisible by all $x_i$ such that $e_i$ is in the
relative interior of $\sigma$, this cocycle is equivalent up to a
coboundary to
 \begin{equation}\label{e:coc}
 {\biggl\{\frac{B}{\prod_{\rho_i\subset\sigma} x_i}
 \biggl(\partial^l_{j_1}-\partial^l_{j_0}\biggr)\biggr\}}_{j_0j_1},
 \end{equation}
 which has an obvious lift by the restriction map
 $H^1(\ps,{\cal T}_\ps)@>>>H^1(X,i^*{\cal T}_\ps)$.
Notice that the lift is independent of the hypersurface $X$, which
is a strong evidence that the non-polynomial deformations are
induced by the deformation of the ambient variety. We should also
point out that the above cocycles are coboundaries in the case
$B\in S_{\beta^\sigma_1}$ is divisible by $x_l$.

\section{Deformations of toric varieties and semiample hypersurfaces.}\label{s:d}

In this section, we  take the open covering of the toric variety
used for the cocycles representing the non-polynomial deformations
and reglue the open sets in a certain way so that  the complex
structure deforms on the toric variety. The construction of a flat
family corresponding to this gluing is automatic. Then,  we also
find a subfamily  which gives a non-polynomial  deformation of
semiample hypersurfaces. We consider complete toric varieties over
complex numbers, but everything holds over an algebraically closed
field.

Let $\ps$ be a complete toric variety with a big and nef divisor
class $[X]\in A_{d-1}(\ps)$.
 Take the open covering
$\{U_{\sigma_j}\}_{j=1}^{n(\sigma)}$, considered in
Section~\ref{s:i} for $\sigma\in\Sigma_X(2)$ with at least one
$e_l$, $l=l_j$ as in (\ref{e:pic}), lying in its relative
interior. It was practical to use this open covering in \cite{m2}
for calculations, but we should note that this covering of the
toric variety $\ps$ is a refinement of a covering by two open sets
$U^l_0$ and $U^l_1$ given by $\prod_{k>j}x_{l_k}\ne0$ and
 $\prod_{k<j}x_{l_k}\ne0$, respectively.
These sets intersect in the open toric subvariety
  $$U^l_0\cap U^l_1=\{x\in\ps:\prod_{\rho_l\ne\rho_i\subset\sigma}x_i\ne0\}.$$

 Next, note that the noncomplete toric variety $U^l_0\cap U^l_1$
  has the following  automorphisms. The monomials
$\prod_{\rho_i\subset\sigma}x_i\prod_{i=1}^n x_i^{\langle
u,e_i\rangle}$ in $S_{\beta^\sigma_1}$, where
$\beta^\sigma_1=\sum_{\rho_i\subset\sigma}\deg(x_i)$, correspond
to $u\in M$ such that $\langle u,e_i\rangle\ge-1$ for
$\rho_i\subset\sigma$ and $\langle u,e_i\rangle\ge0$ for
$\rho_i\not\subset\sigma$.   Then the monomials $x_l\prod_{i=1}^n
x_i^{\langle u,e_i\rangle}$ have the same degree as $x_l$ and do
not have poles in the intersection $U^l_0\cap U^l_1$. The
1-parameter subgroups of automorphisms ${\rm Aut}(U^l_0\cap
U^l_1)$ are induced by
\begin{equation}\label{e:tr}
y_\lambda(x_1,\dots,x_l,\dots,x_n)=(x_1,\dots,x_l+\lambda
x_l\prod_{i=1}^n x_i^{\langle u,e_i\rangle},\dots,x_n).
\end{equation}
 For those  $u$ above, which satisfy $\langle u,e_l\rangle=-1$, the
 lattice point $-u$ is called a {\it root} of the noncomplete toric
 variety $U^l_0\cap U^l_1$ (this is similar to the definitions in
 \cite{c1} and \cite{de}).

 To deform the complex structure on the toric variety $\ps$, we
 reglue the open subsets $U^l_0$ and $U^l_1$  of $\ps$ along $U^l_0\cap U^l_1$
by the open immersions:
 $$U^l_0\cap U^l_1\hookrightarrow U^l_k,\quad k=1,2,$$
which  send a point in
 $U^l_0\cap U^l_1$
 represented by $(x_1,\dots,x_l,\dots,x_n)$ to the corresponding points in
$U^l_k$  represented by $$(x_1,\dots,x_l+(-1)^k\lambda
x_l\prod_{i=1}^n x_i^{\langle u,e_i\rangle},\dots,x_n).$$
  All this is well defined because there are natural commutative
  diagrams:
$$\begin{CD} \{x\in\AA^n\setminus{\bf
V}(B(\Sigma)):\prod_{k<j}x_{l_k}\ne0\}@>>>U^l_0\\
 @AAA @AAA\\
 \{x\in\AA^n\setminus{\bf V}(B(\Sigma)):\prod_{\rho_l\ne\rho_i\subset\sigma}x_i\ne0\}@>>>
U^l_0\cap U^l_1 \\ @VVV @VVV\\ \{x\in\AA^n\setminus{\bf
V}(B(\Sigma)):\prod_{k>j}x_{l_k}\ne0\}@>>>U^l_1,
\end{CD}
$$
 where the horizontal arrows are toric morphisms (see \cite[page~27]{c1})
 and the vertical arrows are the open immersions.
 If $-u$ is {\it not} a root of $U^l_0\cap U^l_1$, it is not
 hard to see that $\prod_{i=1}^n x_i^{\langle u,e_i\rangle}$ will have
 no poles on one of the open sets $U^l_0$ or $U^l_1$, whence (\ref{e:tr})
 induces an automorphism on that set, and the above gluing
 produces the same variety $\ps$.
If $-u$ is the root, the glued set is
 not a toric variety any more, but it still has a Zariski open covering
 by noncomplete toric varieties, in fact, the same covering as the
 initial complete toric variety had. Such varieties fit into the
 category of toroidal varieties (see \cite{d}).

The description of the deformation of complex structure we
presented is analogous to  Spencer's idea on  the complex
deformation  in \cite{ks}. However, algebraic geometry defines a
({\it global}) {\it deformation} of a scheme as a flat family with
one of the fibers isomorphic to the scheme. The construction of
such a flat family from the gluing condition is straightforward.
In our situation, consider two sets $U^l_0\times{\AA}^1$ and
$U^l_1\times{\AA}^1$, which we glue along $U^l_0\cap
U^l_1\times{\AA}^1$ to form a reduced scheme
 ${\cal P}$
by the identification:
\begin{multline}\label{e:id}
 U^l_0\cap U^l_1\times{\AA}^1\hookrightarrow U^l_k\times{\AA}^1,\quad k=1,2,\\
(x_1,\dots,x_l,\dots,x_n,\lambda)\mapsto
(x_1,\dots,x_l+(-1)^k\lambda x_l\prod_{i=1}^n x_i^{\langle
u,e_i\rangle},\dots,x_n,\lambda),
\end{multline}
 where $\lambda$ is an affine coordinate on $\AA^1$.
The obvious projection onto the last component gives us a flat
family ${\cal P}@>>>\AA^1$ whose fiber over the point $\lambda=0$
is precisely the original complete toric variety $\ps$.

 To
describe the deformations of big and nef hypersurfaces
$X\subset\ps$ induced by the above deformations of the ambient
toric variety, we will construct three families in $U^l_0\cap
U^l_1\times{\AA}^1$, $U^l_0\times{\AA}^1$ and
$U^l_1\times{\AA}^1$, which will be patched together by the
identification (\ref{e:id}). For simplicity, we only restrict to
nontrivial deformations  assuming that $-u$ is the root of
$U^l_0\cap U^l_1$.
 Let
$X$ be defined by the polynomial
 $$f(x)=\sum_{m\in\Delta}a_m \prod_{i=1}^n
x_i^{b_i+\langle m,e_i\rangle}$$ in $S_\beta$, where
$\Delta=\Delta_D$ is the polytope in $M$ associated
 to the torus invariant divisor
$D=\sum_{i=1}^n b_iD_i$ and given by the conditions $b_i+\langle
m,e_i\rangle\ge0$. Consider the hypersurface in $U^l_0\cap
U^l_1\times{\AA}^1$ defined by  the equation:
$$
f_\lambda(x):=\sum_{m\in\Delta}a_m \prod_{i=1}^n x_i^{b_i+\langle
m,e_i\rangle}\biggl(1+c_m\lambda\prod_{i=1}^n
 x_i^{\langle u,e_i\rangle}\biggr)^{b_l+\langle m,e_l\rangle}
 =0,
 $$
 where
 $$c_m= \left\{
\begin{matrix}
 1 & \mbox{ if } & b_{l_0}+\langle m,e_{l_0}\rangle+\langle u,e_{l_0}\rangle(b_l+\langle
 m,e_l\rangle)>0
\\
 0 & \mbox{ if } & b_{l_0}+\langle m,e_{l_0}\rangle+\langle u,e_{l_0}\rangle(b_l+\langle
 m,e_l\rangle)=0
 \\
  -1 & \mbox{ if } & b_{l_0}+\langle m,e_{l_0}\rangle+\langle u,e_{l_0}\rangle(b_l+\langle
 m,e_l\rangle)<0.
\end{matrix}\right.$$
Then denote by $f^k_\lambda(x)$ the rational function:
 $$ \sum_{m\in\Delta}a_m \prod_{i=1}^n
x_i^{b_i+\langle m,e_i\rangle}\biggl(1+
((-1)^k+c_m)\lambda\prod_{i=1}^n
 x_i^{\langle u,e_i\rangle}\biggr)^{b_l+\langle m,e_l\rangle},$$
  obtained from $f_\lambda(x)$ by the transformation
(\ref{e:id}).  We claim that this function does not have poles
 on $U^l_k\times{\AA}^1$, and, therefore, defines a hypersurface there.
  Indeed, the only poles can occur along the divisors $D_i$
 for $\rho_i\subset\sigma$. Since
  $\langle u,e_l\rangle=-1$,  there are obviously no poles along $D_l$.
Further, notice that $b_i= \langle m_\sigma,e_{i}\rangle$ for
$\rho_i\subset\sigma$ and some  $m_\sigma\in M$, by the
construction of the fan $\Sigma_X$ (see \cite[Section~1]{m1}).
Hence,  $b_{l_0}+\langle m,e_{l_0}\rangle+\langle
u,e_{l_0}\rangle(b_l+\langle
 m,e_l\rangle)<0$ if and only if
 $b_{l_i}+\langle m,e_{l_i}\rangle+\langle
u,e_{l_i}\rangle(b_l+\langle
 m,e_l\rangle)\ge0$ for any $0\le i< j$,
  because
 $b_{l}+\langle m,e_{l}\rangle+\langle u,e_{l}\rangle(b_l+\langle m,e_l\rangle)=0$.
If $f^0_\lambda(x)$ had a pole along $D_{l_i}$ for $0\le i\le j$,
then  the smallest degree in $x_{l_i}$ would be $b_{l_i}+\langle
m,e_{l_i}\rangle+\langle u,e_{l_i}\rangle(b_l+\langle
 m,e_l\rangle)<0$ for some $m\in\Delta$. However, $c_m=-1$ in this case by the above, which means
 that  $f^0_\lambda(x)$ has no poles on $U^l_0\times{\AA}^1$. A
 similar argument shows that $f^1_\lambda(x)$ defines a hypersurface in $U^l_1\times{\AA}^1$.

  By the identification (\ref{e:id}), the families given by $f^k_\lambda(x)$
 in $U^l_k\times{\AA}^1$, for
$k=1,2$, fit together along
 $f_\lambda(x)=0$ in  $U^l_0\cap
U^l_1\times{\AA}^1$ to form a non-polynomial deformation  ${\cal
X}@>>>\AA^1$ of the semiample hypersurface $X\subset\ps$. This
family is flat because algebraic families of divisors are flat
(see \cite[Chapter III, Example~9.8.5]{h}).

To confirm our construction we want to compare it with the
non-polynomial deformation of a Calabi-Yau hypersurface in the
weighted projective space $\PP(1,1,2,2,2)$ described in
\cite{kklm}.

\begin{ex} Let $Y$ be the quasismooth Calabi-Yau hypersurface in
$\PP(1,1,2,2,2)$,
with homogeneous coordinates $z_1,z_2,z_3,z_4,z_5$,
 defined  by the equation
 $$z_1^ 8 + z_2 ^8 + z_3^ 4
+ z_4^ 4 + z_5^ 4 = 0.$$ The weighted projective space
$\PP(1,1,2,2,2)$, whose fan  has the  following integral
generators of
 the 1-dimensional cones
\begin{equation}\label{e:fan}
\{(-1,-2,-2,-2),(1,0,0,0),(0,1,0,0),(0,0,1,0),(0,0,0,1)\},
\end{equation}
 is singular along $\PP^2$ given by $z_1=z_2=0$. This singularity
corresponds to the fact that the 2-dimensional cone $\sigma$
generated by the first two lattice points $(-1,-2,-2,-2)$ and
$(1,0,0,0)$ contains another lattice point not integrally
generated by these generators:
 $$(0,-1,-1,-1)=\frac{1}{2}(-1,-2,-2,-2)+\frac{1}{2}(1,0,0,0).$$
The crepant desingularization of $\PP(1,1,2,2,2)$ is the toric
variety
 $\ps$ whose fan is obtained from the fan determined by (\ref{e:fan}) by
 inserting the ray through $(0,-1,-1,-1)$.
If $x_1,x_2,x_3,x_4,x_5,x_6$ denote the homogeneous coordinates of
$\ps$, corresponding to the five points in (\ref{e:fan}) and the
sixth  $(0,-1,-1,-1)$, then the crepant resolution $X$ of $Y$ is
defined by the polynomial $$f(x)=x_1^ 8x_6^4 + x_2 ^8x_6^4 +
x_3^4+ x_4^ 4 + x_5^ 4.$$
 To describe the non-polynomial deformation of $X$ first embed $\PP(1,1,2,2,2)$
 into $\PP^5$ by the degree 2 linear system:
 $$(z_1,z_2,z_3,z_4,z_5)\mapsto(z_1^2,z_2^2,z_1z_2,z_3,z_4,z_5).$$
  Then the
 image of $\PP(1,1,2,2,2)$ is the quadric $y_0y_1=y_2^2$, while
 the image of $Y$ is the intersection of this quadric with the
 quartic
 \begin{equation}\label{e:quar}
 y_0^4 + y_1 ^4 + y_3^ 4
+ y_4^ 4 + y_5^ 4 = 0. \end{equation}
 The complex deformation of $X$, as stated in \cite{kklm},
 is the blow-up of the complete intersections of the quartic (\ref{e:quar}) and a
 family  of rank 4 quadrics, like $y_0y_1=y_2(y_2+2\lambda y_3)$,
 in $\PP^5$ along $y_0=y_2=0$.
We can describe this blow-up as the complete intersection
\begin{equation}\label{e:int}
 y_1z_1=(y_2+\lambda y_3)z_2,\quad y_0z_2=y_2z_1
\end{equation}
in $\PP^5\times\PP^1$, where $z_1,z_2$ are the coordinates on the
$\PP^1$ factor. The toric variety $\ps$, which is
 the blow-up of the  quadric $y_0y_1=y_2^2$ along
$y_0=y_2=0$, is isomorphic to the complete intersection
$y_1z_1=y_2z_2$ and $y_0z_2=y_2z_1$ under the map
\begin{equation}\label{e:map}
(x_1,x_2,x_3,x_4,x_5,x_6)\mapsto(x_1^2x_6,x_2^2x_6,x_1x_2x_6,x_3,x_4,x_5)\times(x_1,x_2).
\end{equation}
Comparing  the above deformation with ours, let $l=6$ and $l_0=1$.
Then notice that the charts $U^l_0$ and $U^l_1$, defined by
$x_2\ne0$ and $x_1\ne0$ in $\ps$, respectively, can be identified
with the open subsets $z_2\ne0$ and $z_1\ne0$ of the complete
intersection (\ref{e:int}) by composing (\ref{e:map}) with the
maps, which send $(y,z)$ to  $(y_0-2\lambda
y_3\frac{z_1}{z_2},y_1,y_2-2\lambda y_3,y_3,y_4,y_5,z_1,z_2)$ and
$(y_0,y_1+2\lambda y_3\frac{z_2}{z_1},y_2,y_3,y_4,y_5,z_1,z_2)$,
respectively. This gives the transition function between $U^l_0$
and $U^l_1$: $$(x_1,\dots,x_5,x_6)\mapsto
(x_1,\dots,x_5,x_6-2\lambda \frac{x_3}{x_1x_2}),$$
 precisely as we had in (\ref{e:id}) with $u=(-1,1,0,0)$.
 Under the above identification, the quartic
(\ref{e:quar}) maps to the hypersurfaces in $U^l_0$ and $U^l_1$
determined by the rational functions $x_1^ 8x_6^4(1+2\lambda
\frac{x_3}{x_1x_2}))^4 + x_2 ^8x_6^4 + x_3^4+ x_4^ 4 + x_5^ 4$ and
$x_1^ 8x_6^4+ x_2 ^8x_6^4(1-2\lambda \frac{x_3}{x_1x_2}))^4  +
x_3^4+ x_4^ 4 + x_5^ 4$, which in this case coincide with our
$f_\lambda^k(x)$.
\end{ex}

\section{Matching the global and infinitesimal deformations.}

Here, we find the correspondence between the deformations
constructed in the previous section and the \v Cech cocycles
representing infinitesimal deformations of toric varieties and
infinitesimal non-polynomial deformations of semiample
hypersurfaces in Section~\ref{s:i}.

To compute the infinitesimal deformation of the toric variety, we
can use the formulas in \cite[\S4.2]{ko}. The open covering of
each fiber in the family ${\cal P}@>>>\AA^1$ consists of the two
open sets  $U^l_0$ and $U^l_1$, and the transition function is
given in (\ref{e:id}). Hence, by (4.10) in \cite{ko}, the
corresponding cocycle is $$\biggl\{2\biggl(\prod_{i=1}^n
x_i^{\langle u,e_i\rangle}\biggr)x_l\frac{\partial}{\partial
x_{l}}\biggr\}_{U^l_0U^l_1}.$$
 Since the identity
$$\frac{x_{l_{j-1}}\partial_{l_{j-1}}}{\mult(\sigma_j)}+
\frac{x_{l_{j+1}}\partial_{l_{j+1}}}{\mult(\sigma_{j+1})}
=\frac{\mult(\sigma_j+\sigma_{j+1})}{\mult(\sigma_j)
\mult(\sigma_{j+1})}x_{l_j}\partial_{l_j}$$ holds (see (7) in
\cite{m2}), we conclude that the above cocycle represents the same
element in $H^1(\ps,{\cal T}_\ps)$ as the $-2\mult(\sigma_j)
\mult(\sigma_{j+1})/\mult(\sigma_j+\sigma_{j+1})$  multiple of
 (\ref{e:coc}) with $B=\prod_{\rho_i\subset\sigma} x_i\prod_{i=1}^n x_i^{\langle
u,e_i\rangle}$ does.

For the hypersurface $X\subset\ps$, the infinitesimal deformation
corresponding to the flat family ${\cal X}@>>>\AA^1$ can be found
using \cite{a} and \cite[Chapter III, Section 9]{h}. We replace
the base $\AA^1$ by $\spe(\CC[\varepsilon])$ (abusing notation,
$\CC[\varepsilon]$ denotes
$\CC[\varepsilon]/\langle\varepsilon^2\rangle$): take the morphism
$$\spe(\CC[\varepsilon])@>>>\spe\CC[\lambda]\cong\AA^1,$$ which
arises from the ring homomorhism that sends $\lambda$ to
$\varepsilon$, then by base extension we obtain another family
$${\cal X}'={\cal X}\times_{\AA^1} \spe(\CC[\varepsilon])$$ flat
over the new base. It is not difficult to see that this family is
 glued from the hypersurfaces
  $$f^k_\varepsilon(x)= \sum_{m\in\Delta}a_m \prod_{i=1}^n
x_i^{b_i+\langle m,e_i\rangle}\biggl(1+\varepsilon
((-1)^k+c_m)({b_l+\langle m,e_l\rangle})\prod_{i=1}^n
 x_i^{\langle u,e_i\rangle}\biggr)=0$$
in $U^l_k\times \spe(\CC[\varepsilon])$  along the hypersurface $$
f_\varepsilon(x)=\sum_{m\in\Delta}a_m \prod_{i=1}^n
x_i^{b_i+\langle m,e_i\rangle}\biggl(1+\varepsilon
c_m({b_l+\langle m,e_l\rangle})\prod_{i=1}^n
 x_i^{\langle u,e_i\rangle}\biggr)=0
 $$
 in
$U^l_0\cap U^l_1\times \spe(\CC[\varepsilon])$ with the
identification:
 $$ U^l_0\cap
U^l_1\times\spe(\CC[\varepsilon])\hookrightarrow
U^l_k\times\spe(\CC[\varepsilon]),
 \quad k=1,2,$$
which is defined by sending a function $a+b\varepsilon$ on the
affine scheme $U_{k}^l\cap{\bf
A}_\tau\times\spe(\CC[\varepsilon])$,
 where ${\bf
A}_\tau=\spe(S_\tau)_0$ is the affine toric variety associated
with a cone $\tau\in\Sigma$ and $S_\tau$ is the localization of
$S$ at $\prod_{\rho_i\not\in\tau} x_i$ (see \cite[Lemma~2.2]{c1}),
to the function
\begin{equation}\label{e:tr1}
a+\biggl(b+(-1)^k\biggl(\prod_{i=1}^n
 x_i^{\langle u,e_i\rangle}\biggr)x_l\partial_l a\biggr)\varepsilon
 \end{equation}
on the affine scheme $ U^l_0\cap
U^l_1\times\spe(\CC[\varepsilon])$.

 To find the \v Cech cocycle corresponding the
infinitesimal deformation we  first build trivializations:
  $$\phi_{(i_0,k_0)}:X\cap U_{i_0}\cap
U_{k_0}^l\cap{\bf A}_\tau\times\spe(\CC[\varepsilon])\cong{\cal
X}'\cap U_{i_0}\cap U_{k_0}^l\cap {\bf
A}_\tau\times\spe(\CC[\varepsilon])$$
  by
mapping a
 function $a+b\varepsilon$ on the target scheme to the function
\begin{equation}\label{e:tri}
 a+\biggl(b-\biggl(\prod_{i=1}^n x_i^{\langle
u,e_i\rangle}\biggr)\sum_{m\in\Delta}a_m((-1)^{k_0}+c_m)
\biggl(x_l\partial_l \prod_{i=1}^nx_i^{b_i+\langle
m,e_i\rangle}\biggr) \frac{\partial_{i_0}}{f_{i_0}}
a\biggr)\varepsilon
\end{equation}
on the first scheme. It is
straightforward to check that this map is a well-defined
isomorphism of the rings of functions, which sends
$f^k_\varepsilon/x^D$ to  $f/x^D$ where
 $x^D\in S_\beta$  is a monomial invertible in $S_\tau$ (see
 \cite[page~317]{bc}). The inverse of $\phi_{(i_0,k_0)}$ is given by
 changing the sign after $b$ in (\ref{e:tri}).
Then, taking into account (\ref{e:tr1}), the automorphisms
 $\phi_{(i_0,k_0)}^{-1}\circ\phi_{(i_1,k_1)}$ on
$$X\cap U_{i_0}\cap U_{k_0}^l\cap U_{i_1}\cap U_{k_1}^l\cap{\bf
A}_\tau\times\spe(\CC[\varepsilon])$$ correspond on the level of
rings  to mapping a function $a+b\varepsilon$ to
$a+(b+\theta_{(i_0,k_0)(i_1,k_1)}a)\varepsilon,$ where
$\theta_{(i_0,k_0)(i_1,k_1)}$ is the polyvector field
\begin{multline*}
\biggl(\prod_{i=1}^n x_i^{\langle
u,e_i\rangle}\biggr)\biggl(((-1)^{k_0}-(-1)^{k_1})x_l\partial_l\\
-\sum_{m\in\Delta}a_m \biggl(x_l\partial_l
\prod_{i=1}^nx_i^{b_i+\langle m,e_i\rangle}\biggr)
\biggl(((-1)^{k_0}+c_m)
\frac{\partial_{i_0}}{f_{i_0}}-((-1)^{k_1}+c_m)\frac{\partial_{i_1}}{f_{i_1}}\biggr)\biggr)
\end{multline*}
The \v Cech cocycle
$\{\theta_{(i_0,k_0)(i_1,k_1)}\}_{(i_0,k_0)(i_1,k_1)}$ represents
the infinitesimal deformation ${\cal X}'$ of the hypersurface $X$
in $H^1(X,{\cal T}_X)$. It is not hard to see that, up to a
polynomial infinitesimal deformation $\gamma_A$ for appropriate
$A\in S_\beta$, the cocycle
$$-2\frac{\mult(\sigma_j)\mult(\sigma_{j+1})}{\mult(\sigma_j+\sigma_{j+1})}
(\gamma^l_B)_{(i_0,{j_0})(i_1,{j_1})}$$ with
$B=\prod_{\rho_i\subset\sigma} x_i\prod_{i=1}^n x_i^{\langle
u,e_i\rangle}$ corresponds to the same element in $H^1(X,{\cal
T}_X)$ as $\{\theta_{(i_0,k_0)(i_1,k_1)}\}_{(i_0,k_0)(i_1,k_1)}$
does.

\end{document}